\documentclass[11pt,a4paper]{article}
\long\def\meta#1{\texttt{#1}}
\long\def\drop#1{}

\usepackage[latin1]{inputenc}
\usepackage{amsmath}
\usepackage{amsthm}
\usepackage{amssymb}
\usepackage{latexsym}
\usepackage{graphicx}
\usepackage{mathrsfs}
\usepackage{geompsfi}
\usepackage{psfrag} 
\usepackage{subfig}
\usepackage{bm}

\newcommand{\F}{\mathscr{F}}

\def\R{\mathbb R}
\def\div{\mathop{\mathrm{div}}}
\let\o\Omega
\let\ds\displaystyle

\let\es\varepsilon

\newcommand{\hf}{\mathscr{H}^1}
\def\pref#1{(\ref{#1})}
\newtheorem{theorem}{Theorem}

\newenvironment{remark}%
  {\par\medbreak\refstepcounter{theorem}%
    \noindent\textbf{Remark~\thetheorem. }}%
  {\par\medskip}

\def\Xint#1{\mathchoice
   {\XXint\displaystyle\textstyle{#1}}%
   {\XXint\textstyle\scriptstyle{#1}}%
   {\XXint\scriptstyle\scriptscriptstyle{#1}}%
   {\XXint\scriptscriptstyle\scriptscriptstyle{#1}}%
   \!\int}
\def\XXint#1#2#3{{\setbox0=\hbox{$#1{#2#3}{\int}$}
     \vcenter{\hbox{$#2#3$}}\kern-.5\wd0}}
\def\dashint{\Xint-}

\begin{document}
\title{Stripe Patterns and the Eikonal Equation}
\author{Mark A. Peletier and Marco Veneroni}
\date{\today}
\maketitle

\begin{abstract}
We study a new formulation for the eikonal equation $|\nabla u| =1$ on a bounded subset of $\R^2$. Considering a field $P$ of orthogonal projections onto $1$-dimensional subspaces, with $\div P \in L^2$, we prove existence and uniqueness for solutions of the equation $P\div P=0$. We give a geometric description, comparable with the classical case, and we prove that such solutions exist only if the domain is a tubular neighbourhood of a regular closed curve. 

This formulation provides a useful approach to the analysis of stripe patterns. It is specifically suited to systems where the physical properties of the pattern are invariant under rotation over 180 degrees, such as systems of block copolymers or liquid crystals.

\end{abstract}
\smallskip

\noindent\textbf{AMS Cl.} 35L65, 35B65.

\section{Introduction}


In this note we study a new formulation of the Eikonal equation which was suggested by an example of stripe patterns arising in \emph{block copolymer melts}. For precise statements of the results, complete proofs and references, we refer to \cite{PeletierVeneroniEikonalTA} and \cite{PeletierVeneroniTR}. 
\medskip

Many pattern-forming systems produce parallel stripes, sometimes straight, sometimes curved. In geology, for instance, `parallel folding' refers to the folding of layers of rock in a manner that preserves the layer thickness but allows for curving of the layers~\cite{BoonBuddHunt07}. In a different context, the convection rolls of the Rayleigh-B\'enard experiment produce striped patterns that may also be either straight or curved (see e.g.~\cite{BodenschatzPeschAhlers00}). 

Block copolymers consist of two covalently bonded, mutually repelling parts (`blocks'). At sufficiently low temperature the repelling forces lead to patterns with a length scale that is related to the length of single polymers. We recently studied the behaviour of an energy that describes such systems, and investigated a limit process in which the stripe width tends to zero~\cite{PeletierVeneroniTR}. In that limit the stripes not only become thin, but also uniform in width, and the stripe pattern comes to resemble the level sets of a solution of the eikonal equation. The rigorous version of this statement, in the form of a Gamma-convergence result, gives rise to a new formulation of the eikonal equation, in which the directionality of the stripes is represented by line fields rather than by vector fields. Before stating this formulation mathematically we first describe it in heuristic terms.

\subsection{The Eikonal equation}
The eikonal equation has its origin in models of wave propagation, where the equation describes the position of a wave front at different times $t$. For a homogeneous and isotropic medium, in which the wave velocity is constant, the equation can be written in the form
\begin{equation}
\label{eq:EEvector}
|\nabla u | = 1.
\end{equation}
The wave front at time $t$ is given by the level set $\{x: u(x) = t\}$, and the function $u$ has the interpretation of the time needed for a wave to arrive at the point $x$. 

A feature of the eikonal equation is that the fronts at different times are parallel, provided no singularities occur. In this sense the equation is a natural candidate for the description of other processes that involve parallellism, such as the stripe-forming systems mentioned above. However, a major difference between the stripe-forming systems and the wave-front model is that the wave front has a natural \emph{directionality} associated with it: of the two directions normal to a front, one is `forward in time' and the other `backward'. This distinction also is visible in the notion of viscosity solution for~\pref{eq:EEvector} (see e.g.~\cite{CrandallEvansLions84}).

\begin{figure}[ht]
\centering
\noindent
\includegraphics[height=2cm]{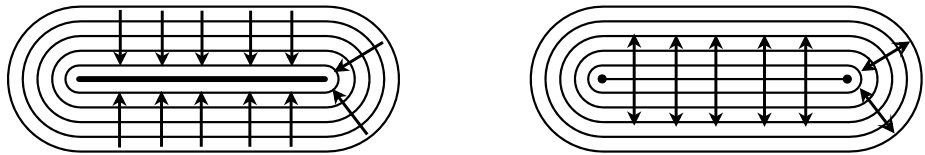}
\caption{Stripe patterns can be represented by vectors (left) or by unoriented line fields (right). Both representations have a vortex singularity at the two ends; but the vector representation also contains a jump singularity along the connecting line. Note that the regularity restrictions of this paper exclude both types of singularity, however.}
\label{fig:stadium}
\end{figure}

The stripe patterns, on the other hand, have no inherent distinction between the two normal directions. As a consequence a vector representation of a stripe pattern may have singularities that have no physical counterpart. Figure~\ref{fig:stadium} (left) shows an example of this.

\subsection{Diblock Copolymers}
An \emph{AB diblock copolymer} is constructed by grafting two polymers together (called the A and B parts). Repelling forces between the two parts lead to phase separation at a scale that is no larger than the length of a single polymer. In this micro-scale separation patterns emerge, and it is exactly this pattern-forming property that makes block copolymers technologically useful~\cite{RuzetteLeibler05}. 

In \cite{PeletierVeneroniTR} we study the formation of stripe-like patterns in a specific two-dimensional system that arises in the modelling of block copolymers. This system is defined by an energy $\mathcal G_\es$ that admits locally minimizing stripe patterns of width $O(\es)$. As $\es\to0$, we show that any sequence $u_\es$ of patterns for which $\mathcal G_\es(u_\es)$ is bounded becomes stripe-like. In addition, the stripes become increasingly straight and uniform in width.

The energy functional, derived in~\cite[Appendix A]{PeletierRoegerTA}, is   
\begin{equation}
\label{eq:functional}
\F_\es(u) = \left\{ \begin{array}{ll} \displaystyle
\es\int_\Omega |\nabla u| + \frac1\es d(u,1-u),
  \qquad&\mbox{ if $u \in K$,}\vspace{0.25cm}\\
\infty & \mbox{ otherwise.} \end{array} \right.
\end{equation}
Here $\Omega$ is an open, connected, and bounded subset of $\R^2$ with $C^2$ boundary, $d$ is the Monge-Kantorovich distance, and
$$		K:=\left\{u\in BV(\Omega;\{0,1\})\, :\ \dashint_\Omega u(x)\, dx= \frac{1}2 
		\ \text{ and }\ u = 0  \text{ on }  \partial \Omega\,\right\}.$$
We introduce a rescaled functional $\mathcal G_\es$ defined by
$$\mathcal G_\es(u) := \frac1{\es^2} \Bigl( \F_\es(u) - |\Omega|\Bigr).$$
The interpretation of the function $u$ and the functional $\F_\es$ are as follows.

\begin{figure}[ht]
\centering
\psfig{figure=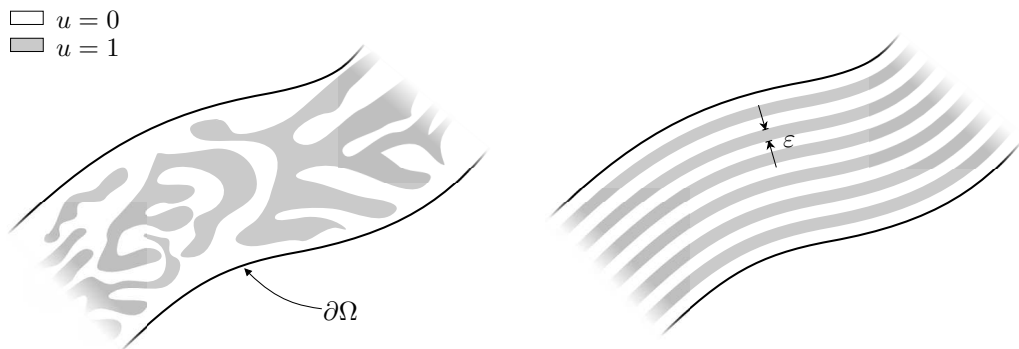,height=5cm}
\caption{A section of a domain $\Omega$ with a general admissible pattern (left) and a stripe-like pattern (right). We prove that in the limit $\es\to0$ all patterns with bounded energy $\mathcal G_\es$ resemble the right-hand picture.}
\label{fig:intro}
\end{figure}

The function $u$ is a characteristic function, whose support corresponds to the region of space occupied by the A part of the diblock copolymer; 
 the complement (the support of $1-u$) corresponds to the B part. The boundary condition $u=0$ in $K$ reflects a repelling force between the boundary of the experimental vessel and the A phase. Figure~\ref{fig:intro} shows two examples of admissible patterns. 

The functional $\F_\es$ contains two terms. The first term penalizes the interface between the A and the B parts, and arises from the repelling force between the two parts; this term favours large-scale separation. In the second term the the Monge-Kantorovich distance~$d$ appears; this term is a measure of the spatial separation of the two sets $\{u=0\}$ and $\{u=1\}$, and favours rapid oscillation. The combination of the two leads to a preferred length scale, which is of order $\es$ in the scaling of~\pref{eq:functional}. 


\subsection{A non-oriented version of Eikonal equation}
A natural mathematical object for the representation of line fields is a \emph{projection}. We define a {projection} to be 
a matrix $P$ that can be written in terms of a unit vector $m$ as $P=m\otimes m$. 
Such a projection matrix has a range and a kernel that are both one-dimensional, and if necessary one can identify a projection $P$ with its range, i.e. with the one-dimensional subspace of $\R^2$ onto which it projects. Note that the independence of the sign of $m$ -- the unsigned nature of a projection -- can be directly recognized in the formula $P=m\otimes m$.

We define $\div P$ as the vector-valued function whose $i$-th component is given by $(\div P)_i:= \sum_{j=1}^2 \partial_{x_j}P_{ij}$. We consider the following problem. Let $\Omega$ be an open subset of $\R^2$.
 Find $P\in L^\infty(\Omega;\R^{2\times 2})$ such that
\begin{subequations}
\label{pb:main}
\begin{eqnarray}
         P^2 = P && \mbox{a.e. in } \o,\label{prpro}\\
         \mbox{rank}(P)=1 && \mbox{a.e. in } \o,\label{prrank}\\
         P \mbox{ is symmetric} && \mbox{a.e. in } \o,\label{prsymm}\\
         \div P \in L^2(\R^2;\R^2) && (\mbox{extended to $0$ outside  }\o),\label{divpint}\\
         P\, \div\! P = 0  && \mbox{a.e. in } \o.\label{pdivpzero}
\end{eqnarray}
\end{subequations}
The first three equations encode the property that $P(x)$ is a projection, in the sense above, at almost every $x$. The sense of property (\ref{divpint}) is that the divergence of $P$ (extended to $0$ outside $\o$), in the sense of distributions in $\R^2$, is an $L^2(\R^2)$ function, which, in particular, implies 
$$        Pn = 0\quad \mbox{in the sense of traces on }\partial\o.$$
The exponent $2$ in~\pref{divpint} is critical in the following sense. Obvious possibilities for singularities in a line field are jump discontinuities (`grain boundaries') and target patterns (see Figure~\ref{fig:types_of_variation}). 

\def\hht{2cm}%
\def\spacing{\hskip1cm}%
\begin{figure}[ht]
\centering
\subfloat[caption][\centering grain boundary]{\includegraphics[height=\hht,clip=on]{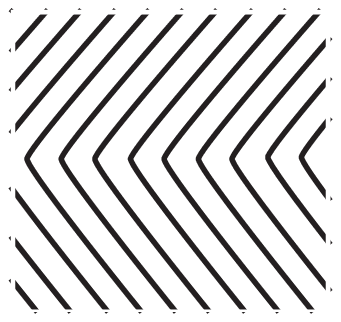}\label{subfig:gb}}
\spacing
\subfloat[caption][\centering target and U-turn patterns]{\includegraphics[height=\hht,clip=on]{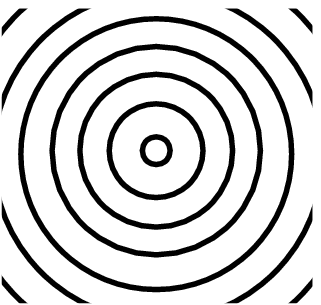}\ \ \includegraphics[height=\hht,clip=on]{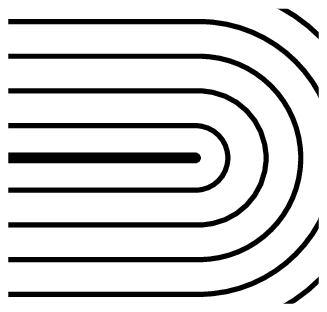}\label{subfig:target}}
\spacing
\subfloat[caption][\centering smooth directional variation]{\includegraphics[height=\hht,clip=on]{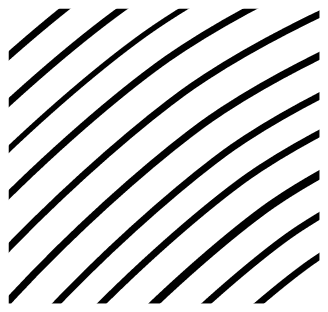}}
\caption{Canonical types of stripe variation in two dimensions. Types~\protect\subref{subfig:gb} and~\protect\subref{subfig:target} are excluded by~\pref{divpint}.
}
\label{fig:types_of_variation}
\end{figure}
At a grain boundary the jump in $P$ causes $\div P$ to have a line singularity, comparable to the one-dimensional Hausdorff measure; condition~\pref{divpint} clearly excludes that possibility. For a target pattern the curvature $\kappa$ of the stripes scales as $1/r$, where $r$ is the distance to the center; then  $\int \kappa^p$ is locally finite for $p<2$, and diverges logarithmically for $p=2$. The cases $p<2$ and $p\geq 2$ therefore distinguish between whether target patterns are admissible ($p<2$) or not.

Given the regularity provided by~\pref{divpint}, the final condition~\pref{pdivpzero} is the eikonal equation itself, as a calculation for a smooth unit-length vector field $m(x)$ shows: 
$$0 = P\div P = m (m\cdot (m\div m + \nabla m \cdot m))
= m\div m + m(m\cdot \nabla m \cdot m)
= m\div m,$$
where the final equality follows from differentiating the identity $|m|^2 =1 $. A solution vector field $m$ therefore is divergence-free, implying that its rotation over 90 degrees is a gradient $\nabla u$; from $|m|=1$ it follows that $|\nabla u|=1$. This little calculation also shows that the 
interpretation of $m$ in $P=m\otimes m$ is that of the stripe direction; $P$ projects along the normal onto the tangent to a stripe.
\medskip

\subsection{Results}
The first result regards the relationship between problem~\pref{pb:main} and the eikonal equation. In the heuristic discussion above, we argued that the field $P$ has to be locally parallel, to be parallel to the boundary of $\Omega$, and to avoid line and vortex singularities. The two first statements are formalized in the following theorem.
\begin{theorem}
\label{th:ee0}
Let $\Omega$ be an open, bounded, and connected subset of $\R^2$ with $C^2$ boundary, and let $P$ be a solution of~\pref{pb:main}.
\begin{enumerate}
\item Let $x_0$ be a Lebesgue point of $P$ in $\Omega$, let $x\in \Omega$, and let $L$ be the line segment connecting $x_0$ with $x$. Assume that $L\subset \Omega$. If $P(x_0)\cdot (x-x_0) =0$, then $P(y)=P(x_0)$ for $\hf$-almost every $y\in L$.
\label{th:ee0:orth}
\item $P\cdot n=0$ a.e. on $\partial \Omega$.
\label{th:ee0:bdry}
\end{enumerate}
These two statements are meaningful since 
\begin{enumerate}
\setcounter{enumi}2
\item $P\in H^1(\Omega;\R^{2\times2})$.
\label{th:ee0:regularity}
\end{enumerate}
\end{theorem}

\bigskip

In the second statement we provide a strong characterization of the \emph{geometry of the domain}~$\Omega$:

\begin{theorem}
\label{th:ee}
Let $\Omega$ be an open, bounded, and connected subset of $\R^2$ with $C^2$ boundary. Then there exists a solution of~\pref{pb:main} if and only if $\Omega$ is a tubular domain. In that case the solution is unique.
\end{theorem}

\begin{figure}[ht]
\centering
\subfloat[caption][\centering The tubular domain $\o$]{\includegraphics[height=\hht,clip=on]{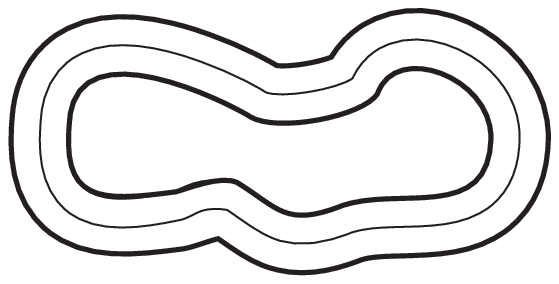}}
\spacing
\subfloat[caption][\centering The projections field $P$]{\includegraphics[height=\hht,clip=on]{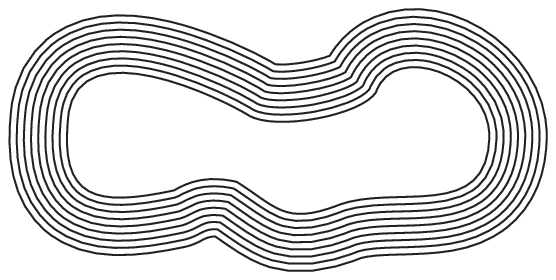}}
\caption{An example of admissible domain, for $P\in \mathcal K_0$.}
\end{figure}
\medskip

A \emph{tubular domain} is a domain in $\R^2$ that can be written as
\[
\Omega = \Gamma + B(0,\delta),
\]
where $\Gamma$ is a closed curve in $\R^2$ with continuous and bounded curvature $\kappa$, $0<\delta< \|\kappa\|^{-1}_\infty$, and $B(0,\delta)$ is the open ball of center 0 and radius $\delta$.

\begin{figure}[ht]
\indent\centering
\psfig{figure=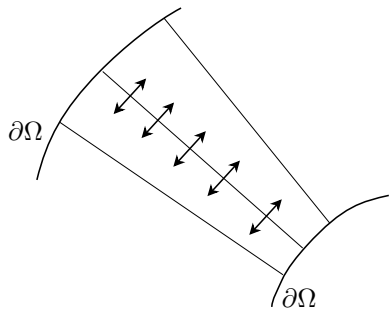,height=4cm}
\caption{If tangent directions propagate normal to themselves (Theorem \ref{th:ee0}.\ref{th:ee0:orth}), and if in addition the boundary is a tangent direction, then the domain is tubular (Theorem~\ref{th:ee}).}
\label{fig:ee}
\end{figure}

The reason why Theorem~\ref{th:ee} is true can heuristically be recognized in a simple picture. Figure~\ref{fig:ee} shows two sections of $\partial \Omega$ with a normal line that connects them. By the first assertion of Theorem~\ref{th:ee0}, the stripe tangents are orthogonal to this normal line; by the second, this normal line is orthogonal to the two boundary segments, implying that the two segments have the same tangent. Therefore the length of the connecting normal line is constant, and as it moves it sweeps out a full tubular neighbourhood. 

The precise relation between the solutions of the non-oriented eikonal equation and the block copolymers energy functionals is the following:
\begin{theorem} The rescaled functional $\mathcal G_\es$ Gamma-converges to the functional 
$$
\mathcal G_0(P):=\left\{\begin{array}{cl}
	\ds \frac 18 \int_\o |\div P(x)|^2dx & \mbox{if } P\in \mathcal{K}_0(\o)\\
	+\infty & \mbox{otherwise}
\end{array} \right.
$$
\end{theorem}
The topology of the Gamma-convergence in this case is the strong topology of measure-function pairs in the sense of Hutchinson \cite{Hutchinson86}. We refer to \cite{PeletierVeneroniEikonalTA, PeletierVeneroniTR} for the details of this proof and an extended discussion.

\subsection{Discussion and future problems}

The work of this paper represents a first step in the analysis of this projection-valued eikonal  equation. While the main results are still lacking in various ways---which we discuss in more detail below---the main point of this paper is to show that this projection-valued formulation is a useful alternative to the usual vector-based formulation. 

To start with, our Theorems~\ref{th:ee0} and~\ref{th:ee} show that solutions of~\pref{pb:main} behave much like we expect from the eikonal equation, in the sense that directional information is preserved in the normal direction. Theorem~\ref{th:ee} makes this property even more explicit, by showing that a full tube, or bunch, of parallel `stripes' can be identified. 

However, it is the differences with the vector-valued eikonal equation that are the most interesting. Figure~\ref{fig:stadium} shows how this formulation can be a better representation of the physical reality than the vector-based form. On the left, the vector field has a jump discontinuity along the center line, while on the right the projection is continuous along that line. Depending on the underlying model, this singularity may have a physical counterpart, or may be a spurious consequence of the vector-based description. For the wave-propagation model the singularity is very real; for striped-pattern systems it typically is not. A projection-valued formulation therefore provides an alternative to the Riemann-surface approach that is sometimes used~\cite{ErcolaniIndikNewellPassot03}. For this distinction to have any consequence, however, solutions with less regularity than the $\div P\in L^2$ of this paper are to be considered. 

\medskip
The proof of the properties stated in Theorem \ref{th:ee0} and Theorem \ref{th:ee} relies on a reduction of the projection-valued formulation to a vector-based formulation which allows us to apply a generalized method of characteristics, introduced in \cite{JabinOttoPerthame02}, to a suitable vector field $m$ satisfying $P=m\otimes m$. This reduction to a vector formulation is achieved by a Lemma by Ball and Zarnescu \cite{BallZarnescu07TA}, which requires $\div P\in L^2$; for less regularity the existence of a lifting may not hold, as the example of the U-turn pattern (Figure~\ref{subfig:target}) shows.  The dependence of the proof on a vector-based representation is awkward in various ways. To start with, the condition $\div P\in L^2$ required for the lifting is much stronger than the conditions that Jabin, Otto, and Perthame require for their results~\cite{JabinOttoPerthame02}. It also has the effect of excluding all singularities, as we already remarked. It would be interesting to prove properties such as those of Theorems~\ref{th:ee0} and~\ref{th:ee} by methods that do not rely on this lifting. 

We would hope that such an intrinsic projection-based proof could also be generalized to the study of target patterns and U-turns, and eventually of grain boundaries. These will require increasingly weak regularity requirements: target patterns may exist for $\div P\in L^p$ with $p<2$, and for a line discontinuity, such as a grain boundary, $\div P$ will be a measure. A natural way to study this kind of patterns could be by exploring different rescalings of the functionals $\F_\es$ defined in \pref{eq:functional}. For example it would be interesting to study the limit of the functionals
$$ \mathcal H_\es(u):=\frac{\F_\es(u) - |\o|}{\es\log \es},$$
which would allow for the formation of target patterns in the limit.

\small
\bibliographystyle{abbrv}
\bibliography{eikon}
\normalsize
\pagebreak

\noindent Mark A. Peletier\smallskip

\noindent Dept. of Mathematics and Computer Science\\
Technische Universiteit Eindhoven\\
PO Box 513\\
5600 MB  Eindhoven\\
The Netherlands\\
e-mail: \meta{m.a.peletier@tue.nl}
\vspace{1cm}

\noindent Marco Veneroni\smallskip

\noindent Technische Universit\"at Dortmund\\
Fakultät für Mathematik, Lehrstuhl I\\
Vogelpothsweg 87\\
44227 Dortmund\\
Germany\\
e-mail: \meta{marco.veneroni@math.uni-dortmund.de}

\end{document}